\documentclass[reqno,a4paper,12pt]{amsart}
\usepackage{amssymb}
\usepackage[english]{babel}
\usepackage[dvips, lmargin=3cm, rmargin=3cm, tmargin=3cm, bmargin=3cm]{geometry}
\usepackage[colorlinks=true,urlcolor=blue,linkcolor=blue]{hyperref}
\usepackage{hyperref}
\usepackage[T1]{fontenc}
\usepackage[english]{babel}

\usepackage{enumitem}
\usepackage{dsfont}
\newtheorem{thm}{Theorem}[section]
\newtheorem{defini}{Definition}[section]
\newtheorem{rem}{Remark}[section]
\newtheorem{lem}{Lemma}[section]
\newtheorem{prop}{Proposition}[section]

\def \R {\mathbb{R} }

\begin{document}
\title[Introduction to Fractional Orlicz-Sobolev spaces]{Introduction to Fractional Orlicz-Sobolev spaces}
\author[Mohammed SRATI, Elhoussine AZROUL and Abdelmoujib Benkirane]
{Mohammed SRATI$^1$\\ Elhoussine AZROUL$^2$\\  Abdelmoujib BENKIRANE$^3$}
\address{[Mohammed SRATI, Elhoussine AZROUL and Abdelmoujib BENKIRANE\newline
 Sidi Mohamed Ben Abdellah
 University,
 Faculty of Sciences Dhar Al Mahraz, Laboratory of Mathematical Analysis and Applications, Fez, Morocco.}
\email{$^1$srati93@gmail.com}
\email{$^2$elhoussine.azroul@gmail.com}
\email{$^3$abd.benkirane@gmail.com}

\maketitle
\textbf{Abstract}
 In this paper, we define the fractional Orlicz-Sobolev spaces,
 and we prove some important results of these spaces. The main result is to show the continuous and compact embedding  for these spaces. As an application, we prove the existence and uniqueness of a solution for a non local problem involving the fractional M-Laplacian operator.
 \section{Introduction}
Let $\Omega$ be an open subset of $\R^N$, and let $s\in (0,1)$. For any $p\in [1,+\infty)$,
the fractional Sobolev spaces are defined as, \\
 $$ W^{s,p}(\Omega)=\Bigg\{u \in L^p(\Omega)  \text{ : }  \dfrac{|u(x)-u(y)|}{|x-y|^{\frac{N}{p}+s}} \in L^p(\Omega \times \Omega)\Bigg\}.$$
 These spaces
have been a classical topic in functional and harmonic analysis all along,
and some important books, such as \cite{36} treat the topic in detail.
On the other hand, fractional spaces, and the corresponding nonlocal equations, are now experiencing impressive applications in different subjects, such as, among others, the thin obstacle problem \cite{31}, finance \cite{17}, phase transitions  \cite{2,10}, stratified
materials \cite{14,15}, crystal dislocation \cite{7}, soft thin films \cite{29}, semipermeable membranes and flame
propagation \cite{11}, conservation laws \cite{8}, ultra-relativistic limits of quantum mechanics \cite{26}, quasi-geostrophic flows \cite{13}, multiple scattering \cite{22}, minimal surfaces \cite{12} , materials science \cite{5}, water
waves \cite{18,19,37},
gradient potential theory \cite{33} and singular set of minima of variational
functionals \cite{32}. Don't panic, instead, see also \cite{35} for further
motivation.\\

 In mathematics and precisely in PDEs, when trying to relax some conditions on the operators (as growth conditions), the problem can not be formulated with classical Lebesgue and Sobolev spaces. Hence, the  adequate functional spaces is the so-called Orlicz spaces (see \cite{3},\cite{4}).\\

 J.F.Bonder and A.M.Salort in \cite{9} proposed a version of the fractional Orlicz-Sobolev spaces, i.e.,
 $$W^s{L_M}(\Omega)=\Bigg\{u\in L_M(\Omega) :  \int_{\Omega} \int_{\Omega} M(\dfrac{ u(x)- u(y)}{|x-y|^s})\dfrac{dxdy}{|x-y|^N}< \infty \Bigg\},$$
 where $M$ is N-function and $L_M$ is the Orlicz space (see section 2).\\
 The authors prove that for any $u \in L_M(\Omega)$ and $0<s<1$, it holds that 
 $$\lim\limits_{s\rightarrow 1}(1-s)\int_{\mathbb{R}^N} \int_{\mathbb{R}^N} M(\dfrac{ u(x)- u(y)}{|x-y|^s})\dfrac{dxdy}{|x-y|^N}= \int_{\R^N}\widetilde{M}(u(x))dx, $$
 where $\widetilde{M}$ is an N-function defined by,
 $$\widetilde{M}(a):= \lim\limits_{s \rightarrow 1 }(1-s)\int_{0}^{1}\int_{\mathbb{S}^{N-s}}M(a|z_N|r^{1-s})dS_z \dfrac{dr}{r}.$$

The previous definition creates problems in the mathematical analysis, more precisely in the calculus and in the embedding results, for example, the Borel measure defined as $d\mu=\frac{dxdy}{|x-y|^N}$ is not finish in the neighbourhood of the origin. To overcome those difficulties, we introduce another definition of the fractional Orlicz-Sobolev space, i.e.,
$$W^s{L_M}(\Omega)=\Bigg\{u\in L_M(\Omega) : \exists \lambda >0 \int_{\Omega} \int_{\Omega} M\left( \dfrac{\lambda( u(x)-u(y))}{|x-y|^sM^{-1}(|x-y|^N)}\right) dxdy< \infty \Bigg\}.$$
 Observe that in the case $M(t)=t^p$, these spaces coincide with the fractional 
 Sobolev space $W^{s,p}(\Omega)$.\\
\hspace*{0.3cm} We begin this paper by showing some natural properties of the space $W^sL_M(\Omega)$, and then we come to the important point of the article, i.e, to study the embedding results of these spaces. 
 We follow the approach of Donaldson and Trudinger in \cite{21} and show the embedding results of the fractional Orlicz-Sobolev spaces $ W^sL_M(\Omega)$ into the Orlicz spaces.\\
\hspace*{0.3cm} This paper is organized as follows:
 We introduce in the second section, some properties on the Orlicz-Sobolev and fractional Sobolev spaces. \\
 The third section is devoted to proving we prove some important results on the fractional Orlicz-Sobolev spaces, and we prove a result of continuous and compact embedding of these spaces into the Orlicz spaces. 
 Finally, we conclude this article by an application of our main results, to show  the existence and uniqueness of solution for a non local problem  involving the fractional M-Laplacian operator.
\section{Some preliminary results}$\label{100}$
First, we briefly recall the definitions and some elementary properties of the Orlicz spaces and Orlicz-Sobolev spaces. We refer the reader to \cite{1,28,34} for further reference and for some of the proofs of the results in this subsection.
\subsection{Orlicz-Sobolev Spaces} 
 We start by recalling the definition of the well-known N-functions.
 
  Let $\Omega$ be an open subset of $\R^N$. Let $M$ :
  $\R^+ \rightarrow \R^+$ be an N-function, i.e., $M$ is continuous, convex, with $M(t) > 0$ for $t >
  0$, $ \frac{M(t)}{t}\rightarrow 0$ as  $t \rightarrow 0$ and $ \frac{M(t)}{t} \rightarrow\infty$ as $t \rightarrow \infty$. Equivalently, $M$ admits the
  representation : $ M(t)=\int_{0}^{t} m(s)ds$ where $m : \R^+ \rightarrow \R^+$ is non-decreasing, right
  continuous, with $m(0) = 0$, $ m(t) > 0$ $\forall t> 0$ and $m(t) \rightarrow \infty $ as $t \rightarrow \infty$.  The conjugate N-function of $M$ is defined by $\overline{M}(t) = \int_{0}^{t} \overline{m}(s)ds$, where $ \overline{m} : \R^+\rightarrow
  \R^+$ is given by $\overline{m}(t) = \sup \left\lbrace s : m(s) \leqslant t\right\rbrace$.\\ Evidently we have
  \begin{equation}\label{1}
st\leqslant M(t)+\overline{M}(s),
  \end{equation}
  which is known Young's inequality. Equality holds in (\ref{1}) if and only if either $t=\overline{m}(s)$ or $s=m(t)$.\\
   We will extend these
  N-functions into even functions on all $\R$. The N-function $M$ is said to satisfy the
  global $\Delta_2$-condition if, for some $k > 0,$
 $$ M(2t) \leqslant kM(t)\text{ , }\forall t \geqslant 0.$$
  When this inequality holds only for $t \geqslant t_0 > 0$, $M$ is said to satisfy the $\Delta_2$-condition
  near infinity.
\\
We call the pair $(M,\Omega)$ is $\Delta$-regular if either :\\
(a) $M$  satisfies a
  global $\Delta_2$-condition, or \\
  (b) $M$ satisfies a $\Delta_2$-condition near infinity and $\Omega$ has finite volume.
\begin{lem}$\label{2.1}$ (cf. \cite{9}).
 Let $M$ be an N-function which satisfies the $\Delta_2$-condition. Then we have,
 \begin{equation}\label{2}
 \overline{M}(m(t))\leqslant (p-1)M(t),
 \end{equation} 
 for some $p>1$. \\
\end{lem}          
Let $\varPhi_1$, $\varPhi_2$ be two N-function. 
       $\varPhi_1$ is stronger (resp essentially stronger) than $\varPhi_2$,  $\varPhi_1\succ\varPhi_2$ (resp $\varPhi_1\succ\succ\varPhi_2$) in symbols, if 
     $$\varPhi_2(x)\leqslant \varPhi_1(a x ), \text{   , } x\geqslant x_0\geqslant 0, $$
     for some (resp for each) $a>0$ and $x_0$ (depending on $a$).
       \begin{rem}
 $\varPhi_1\succ\succ\varPhi_2$  is equivalent to the condition \\
                $$\lim_{x\rightarrow \infty}\dfrac{\varPhi_2(\lambda x)}{\varPhi_1(x)}=0,$$
                for all $\lambda>0$. 
                   \end{rem}
   Let $\Omega$ be an open subset of $\R^N$. The Orlicz class $K_M (\Omega)$ (resp. the Orlicz space
   $L_M(\Omega)$) is defined as the set of (equivalence classes of) real-valued measurable
   functions $u$ on $\Omega$ such that  
\begin{equation}\label{3}
\int_{\Omega} M(u(x))dx <\infty \hspace*{0.5cm} \text{ (resp. } \int_{\Omega} M(\lambda u(x))dx< \infty \text{ for some } \lambda >0 ).
\end{equation}
$L_M(\Omega)$ is a Banach space under the Lexumburg norm
\begin{equation}\label{4}
 ||u||_M=\inf \Bigg\{\lambda>0  : \int_{\Omega}M( \dfrac{u(x)}{\lambda})dx \leqslant 1\Bigg\},
\end{equation}
and $K_M(\Omega)$ is a convex subset of $L_M(\Omega)$. The closure in $L_M(\Omega)$ of the set of
bounded measurable functions on $\Omega$ with compact support in $\overline{\Omega}$ is denoted by $E_M(\Omega)$.\\
The equality $E_M(\Omega)=L_M(\Omega)$ holds if and only if $(M,\Omega)$ is $\Delta$-regular.
 \begin{thm}$\label{2.1,}$ [cf. \cite{1}]
          Let $\Omega$ be an open subset of $\R^N$ which has a finite volume, and suppose  $M,B$ two $N$-function such that $B \prec\prec M$.Then any bounded subset $S$ of $L_{M}(\Omega)$ which is precompact  in $L^1(\Omega)$, is also precompact in $L_B(\Omega)$.\\
  \end{thm}        
      Let $J$ be a nonnegative, real-valued function belonging to $C_0^{\infty}(\R^N)$ and having the properties :
      \begin{itemize}

     \item    $J(x)=0$ if $|x|\geqslant 1$ and 
     \item    $\int_{\R^N}J(x)dx=1$.
        \end{itemize}
       If $\varepsilon>0$, the function 
       $J_{\varepsilon}(x)=\dfrac{1}{\varepsilon^N}J(\dfrac{x}{\varepsilon})$ is nonegative, belongs to $C_0^{\infty}(\R^N)$, and satisfies
       \begin{itemize}

        \item  $J_{\varepsilon}(x)=0$ if $|x|\geqslant 1$, 
        \item  $\displaystyle\int_{\R^N}J_{\varepsilon}(x)dx=1$.
               \end{itemize}
        $J_\varepsilon$ is called a mollifier. We pose
$$
         u_{\varepsilon}(x)=J_{\varepsilon}\ast u(x)=\int_{\R^N}J_{\varepsilon}(x-y)u(y)dy.
$$
      \begin{lem}(cf. \cite{1})\\
       $\bullet$ If $u\in L^1_{loc}(\overline{\Omega})$ then $u_{\varepsilon} \in C^{\infty}(\R^N)$.\\
       $\bullet$ If also supp$(u)$ is compact, then $u_{\varepsilon} \in C_0^{\infty}(\Omega)$, for all $\varepsilon< dist(supp (u), \partial \Omega)$.\\
      
       \end{lem}
      \begin{thm}\label{2.2.}(cf. \cite{1})
      \begin{enumerate}
      \item $C_0(\Omega)$ is dense in $E_M(\Omega)$.
      \item $E_M(\Omega)$ is separable.
      \item $C_0^{\infty}(\Omega)$ is dense in $E_M(\Omega)$.
       \item For each $u \in E_M(\Omega)$, we have $\lim\limits_{\varepsilon \rightarrow 0^{+}}u_{\varepsilon}=u$ in $E_M(\Omega)$.
      \end{enumerate}
      \end{thm}

\subsection{Fractional Sobolev spaces} 
This subsection is devoted to the definition of the fractional Sobolev spaces, and we recall some result of continuous and compact embedding in fractional Sobolev spaces. We refer the reader to \cite{20,24} for further reference and for some of the proofs of these results.

We start by fixing the fractional exponent $s\in(0,1)$. For any $p\in [1,\infty)$, we define the fractional Sobolev space $W^{s,p}(\Omega)$ as follows,
$$ W^{s,p}(\Omega)=\Bigg\{u \in L^p(\Omega)  \text{ : }  \dfrac{|u(x)-u(y)|}{|x-y|^{\frac{N}{p}+s}} \in L^p(\Omega \times \Omega)\Bigg\};$$
i.e, an intermediary Banach space between, endowed
with the natural norm
$$||u||_{s,p}=\Bigg(\int_{\Omega}|u|^pdx+\int_{\Omega}\int_{\Omega} \dfrac{|u(x)-u(y)|^p}{|x-y|^{sp+N}}dxdy\Bigg)^{\frac{1}{p}}.$$
\begin{thm}$\label{2.3}$(cf. \cite{24}).
 Let $s\in (0,1)$ and let $p\in [1,+\infty)$ such that $sp<N$. Let $\Omega$  be an open
  subset of  $\R^N$  with $C^{0,1}$-regularity and bounded boundary. So there exists a constant $C=C(N,s,p,\Omega)$ such that, for all  $f \in W^{s,p}(\Omega)$ we have  \\
   $$||f||_{L^q(\Omega)} \leqslant C||f||_{W^{s,p}(\Omega)} \text{    } \forall q \in [p,p^*],$$
   i.e,
   $$ W^{s,p}(\Omega)  \hookrightarrow L^q(\Omega)  \text{    } \forall q \in[p,p^*],$$
     where $p^*=\frac{Np}{N-sp}$.\\
     If, in addition, $\Omega$ is bounded, then the space $W^{s,p}(\Omega)$ is continuously
     embedded in $L^q(\Omega)$ for any $q\in [1,p^*]$.
\end{thm}
\begin{thm}$\label{2.4}$(cf. \cite{20}).
 Let $s\in (0,1)$ and let $p\in [1,+\infty)$ such that $sp<N$. Let $\Omega$  be a bounded open
  subset of  $\R^N$   with $C^{0,1}$-regularity and bounded boundary. Then the embedding 
 $$ W^{s,p}(\Omega)  \hookrightarrow L^q(\Omega)  \text{    } \forall q \in[1,p^*),$$
 is compact.
\end{thm}
\section{Main results}
\subsection{Fractional Orlicz-Sobolev spaces} 
Now, we define the fractional Orlicz-Sobolev spaces, and we will present some important results on these spaces.
\begin{defini}
Let $M$ be an N-function. For a given domain $\Omega$ in $\R^N$ and $0<s<1$, we define  the fractional Orlicz-Sobolev space $W^sL_M(\Omega)$ as follows,
\small{\begin{equation}\label{5}
W^s{L_M}(\Omega)=\Bigg\{u\in L_M(\Omega) : \exists \lambda >0 / \int_{\Omega} \int_{\Omega} M\left( \dfrac{\lambda( u(x)- u(y))}{|x-y|^sM^{-1}(|x-y|^N)}\right) dxdy< \infty \Bigg\}.
\end{equation}
}
This space is equipped with the norm,
\begin{equation}\label{6}
||u||_{s,M}=||u||_{M}+[u]_{s,M},
\end{equation}

where $[.]_{s,M}$ is the Gagliardo seminorm, defined by 

\begin{equation}\label{7}
[u]_{s,M}=\inf \Bigg\{\lambda >0 :  \int_{\Omega} \int_{\Omega} M\left( \dfrac{u(x)- u(y)}{\lambda|x-y|^sM^{-1}(|x-y|^N)}\right) dxdy\leqslant 1 \Bigg\}.
\end{equation}
\end{defini}
\begin{defini}
Let $M$ be an N-function. For a given domain $\Omega$ in $\R^N$ and $0<s<1$, We define, the space $W^sE_M(\Omega)$ as follows,
\begin{equation}\label{8}
W^s{E_M}(\Omega)=\left\lbrace u\in E_M(\Omega) : \dfrac{|u(x)-u(y)|}{|x-y|^sM^{-1}(|x-y|^N)}\in E_M(\Omega \times \Omega) \right\rbrace.
\end{equation}
\end{defini}
\begin{rem} \text{  }\\
$\bullet$ $W^s{E_M}(\Omega)\subset W^s{L_M}(\Omega)$.\\
 $\bullet$ $W^s{E_M}(\Omega)$ coincides with  $W^s{L_M}(\Omega)$ if and only if  $(M,\Omega)$ is $\Delta$-regular.\\ 
 $\bullet$ If $1<p<\infty$ and $M_p(t)=t^p$, then $W^sL_{M_p}(\Omega)=W^sE_{M_p}(\Omega)=W^{s,p}(\Omega)$.
 \end{rem}
 Many properties of fractional Orlicz-Sobolev spaces are obtained by very straight-forward generalization  of the proofs of the same properties for ordinary fractional Sobolev spaces and Orlicz-Sobolev spaces.

\begin{thm}
Let $\Omega$ be an open subset of $\R^N$, and let $s\in (0,1)$. The space $W^sL_M(\Omega)$ is a Banach space with respect to the norm $(\ref{6})$, and a reflexive  (resp. separable) space if and only if $(M,\Omega)$ is $\Delta$-regular (resp. $(M,\Omega)$ and $(\overline{M},\Omega)$ are $\Delta$-regular).
\end{thm}
 \noindent \textbf{Proof}. 
 Let $\left\lbrace u_n\right\rbrace $ be a Cauchy sequence for the norm $||.||_{s,M}$. In particular, $\left\lbrace u_n\right\rbrace $
 is a Cauchy sequence in $L_M(\Omega)$. It converges to a function $u\in L_M(\Omega)$. Moreover, the
 sequence $\left\lbrace v_n\right\rbrace $ defined as,
 $$v_n(x,y)=\dfrac{|u_n(x)-u_n(y)|}{|x-y|^sM^{-1}(|x-y|^N)},$$
 is a Cauchy sequence in $L_M(\Omega\times \Omega)$. It therefore also converges to an element of  $L_M(\Omega\times \Omega)$.
 Let us extract a subsequence $\left\lbrace u_{\sigma(n)}\right\rbrace $ of $\left\lbrace u_n\right\rbrace $ that converges almost everywhere to $u$. We note that $v_{\sigma(n)}(x, y)$ converges, for almost every pair $(x, y)$, to 
 $$v(x,y)=\dfrac{|u(x)-u(y)|}{|x-y|^sM^{-1}(|x-y|^N)}.$$
 Applying Fatou's lemma, we obtain, for some $\lambda$ (note that $\lambda$ exists since $\left\lbrace u_{\sigma(n)}\right\rbrace  \subset W^sL_M(\Omega)$ ),
{\small$$ \int_{\Omega} \int_{\Omega}  M\left( \dfrac{\lambda (u(x)- u(y))}{|x-y|^sM^{-1}(|x-y|^N)}\right) dxdy\leqslant \liminf_{n\rightarrow \infty} \int_{\Omega} \int_{\Omega} M\left( \dfrac{\lambda( u_{\varphi(n)}(x)- u_{\varphi(n)}(y))}{|x-y|^sM^{-1}(|x-y|^N)}\right)  dxdy<\infty.$$}
Hence $u\in W^s{L_M}(\Omega)$.\\
On the other hand, since $\left\lbrace v_n \right\rbrace$ converges in $L_M(\Omega\times \Omega)$, then by dominated convergence theorem, there exists a subsequence $\left\lbrace v_{\sigma(n)}\right\rbrace $ and a function $ h$ in $L_M(\Omega\times \Omega)$ such that 
$$|v_{\sigma(n)}(x,y)|\leqslant |h(x,y)| \text{ for almost every pair (x,y),} $$
and we have 
$$v_{\sigma(n)}(x,y)\longrightarrow v(x,y)  \text{ for almost every pair (x,y),} $$
this implies by dominated convergence theorem that, $$[u_n-u]_{s,M}\longrightarrow 0.$$
Finally $u_n\rightarrow u$ in $W^sL_M(\Omega)$.\\
\hspace*{0.3cm} To establish the reflexivity and separation of the fractional  Orlicz-Sobolev  spaces, we define the operator T : $W^sL_M(\Omega)\rightarrow L_{M}(\Omega) \times  L_M(\Omega \times \Omega)$ by \\
 $$T(u)=\left( u(x), \dfrac{|u(x)-u(y)|}{|x-y|^sM^{-1}(|x-y|^N)}\right) .$$
 Clearly, T is an isometry. Since $L_M(\Omega)$ is a reflexive and separable space, then $W^sL_M(\Omega)$ is also a reflexive and separable space.\\
 \hspace*{15cm $\Box$ }\\
\hspace*{0.3cm} We are now a position to construct the Orlicz norm corresponding to the fractional Orlicz-Sobolev spaces, and show that it is equivalent to the Lexumburg norm.\\
 The Orlicz norm in Orlicz space is defined by,
 $$||u||_{(M)}=\sup_{\rho(v,\overline{M})\leqslant 1}\left| \int_{\Omega} u(x)v(x)dx\right| ,$$
  where
  $$\rho(v,\overline{M})=\int_{\Omega}\overline{M}(v(x))dx.$$
  By \cite{28}, the expression $||u||_{(M)}$ is a norm in $L_M(\Omega)$ which is equivalent to  $||u||_M$.

 \begin{prop}
  We pose \\
\begin{equation} \label{9}
||u||_{(s,M)}=||u||_{(M)}+[u]_{(s,M)},
\end{equation}
 where
\begin{equation}\label{10}
[u]_{(s,M)}=\sup_{\widetilde{\rho}(v,\overline{M})\leqslant 1 }\left| \int_{\Omega}\int_{\Omega}\dfrac{u(x)-u(y)}{|x-y|^sM^{-1}(|x-y|^N)}v(x,y)dxdy\right|,
\end{equation}
  with 
  $$\widetilde{\rho}(v,\overline{M})=\int_{\Omega}\int_{\Omega}\overline{M}(v(x,y))dxdy.$$
  Then $||.||_{(s,M)}$ is a norm in $W^sL_M(\Omega)$ which is equivalent to $||.||_{s,M}$.
 \end{prop}
  \noindent \textbf{Proof}. 
  Clearly $[.]_{(s,M)}$ is a seminorm, then $||.||_{(s,M)}$ is a norm.\\
 Let  $u \in W^s{L_M}(\Omega)$. By  Young's inequality we have, \\
{\small $$
  \begin{aligned}
  \left[ u\right]_{(s,M)} &=\sup_{\widetilde{\rho}(v,\overline{M}) \leqslant 1 }\left| \int_{\Omega}\int_{\Omega}\dfrac{u(x)-u(y)}{|x-y|^sM^{-1}(|x-y|^N)}v(x,y)dxdy\right|\\
  &\leqslant\sup_{\widetilde{\rho}(v,\overline{M})\leqslant 1 }\left| \int_{\Omega}\int_{\Omega}\left[ M\left( \dfrac{u(x)-u(y)}{|x-y|^sM^{-1}(|x-y|^N)}\right) +\overline{M}\left( v(x,y)\right) \right] dxdy\right|\\
  &\leqslant  \sup_{\widetilde{\rho}(v,\overline{M})\leqslant 1 } \int_{\Omega}\int_{\Omega} M\left( \dfrac{u(x)-u(y)}{|x-y|^sM^{-1}(|x-y|^N)}\right) dxdy+ \sup_{\widetilde{\rho}(v,\overline{M})\leqslant 1 }\int_{\Omega}\int_{\Omega} \overline{M}(v(x,y))dxdy\\
  &\leqslant \int_{\Omega}\int_{\Omega} M\left( \dfrac{u(x)-u(y)}{|x-y|^sM^{-1}(|x-y|^N)}\right)dxdy+1\\
  &:=\phi_{s,M}(u)+1.
  \end{aligned}
  $$}
  Then, we get, 
\begin{equation}\label{11}
\left[ \frac{u}{[u]_{s,M}}\right]_{(s,M)}\leqslant \phi_{s,M}\left( \frac{u}{[u]_{s,M}}\right) +1 \leqslant 2,
\end{equation}
this implies that,  $$[u]_{(s,M)}\leqslant 2[u]_{s,M}.$$\\
  On the other hand, we pose 
  $$ \varphi(x,y)=\dfrac{|u(x)-u(y)|}{|x-y|^sM^{-1}(|x-y|^N)}\in L_M(\Omega\times \Omega).$$ \\
   Then, we have,
   $$\int_{\Omega} \int_{\Omega} M\left( \dfrac{\varphi(x,y)}{||\varphi(x,y)||_{(M)}}\right)dxdy \leqslant 1,$$
   but \\
   $$
   \begin{aligned}
   ||\varphi(x,y)||_{(M)}&=\sup_{\rho(v,\overline{M})\leqslant 1}\left| \int_{\Omega}\int_{\Omega}\varphi(x,y)v(x,y)dxdy\right| \\
   &= \sup_{\widetilde{\rho}(v,\overline{M})\leqslant 1 }\left| \int_{\Omega}\int_{\Omega}\dfrac{u(x)-u(y)}{|x-y|^sM^{-1}(|x-y|^N)}v(x,y)dxdy\right|\\
   &=[u]_{(s,M)}.
   \end{aligned}
   $$
  Finally
   \begin{equation}\label{12}
   \phi_{s,M}\left( \dfrac{u}{[u]_{(s,M)}}\right) =\int_{\Omega} \int_{\Omega} M\left( \dfrac{\varphi(x,y)}{||\varphi(x,y)||_{(M)}}\right)dxdy \leqslant 1,
   \end{equation}
   this implies that  $
[u]_{s,M}\leqslant [u]_{(s,M)}
$.\\
 \hspace*{15cm $\Box$}
\subsection{Approximation theorem and generalized Poincar\'e inequality} 
As in the classic case with $s$ being an integer, any function in the
fractional Orlicz-Sobolev spaces can be approximated by a sequence of
smooth functions with compact support.\\

\begin{thm}\label{3.2.}
 $C^{\infty}_0(\R^N)$ is dense in $W^s{E_M}(\R^N)$.
\end{thm}
  \noindent \textbf{Proof}.  
 Let $u\in W^s{E_M}(\R^N)$, and $\varepsilon>0$. Then  by theorem \ref{2.2.},  $u_{\varepsilon}\in C^{\infty}_0(\R^N) $ and,\\
\begin{equation}
 \lim\limits_{\varepsilon\rightarrow 0^{+}}||u_{\varepsilon}-u||_{(M)}=0.
\end{equation}
So just show that,
  $$  \lim\limits_{\varepsilon\rightarrow 0^{+}}[u_{\varepsilon}-u]_{(s,M)}=0.$$
  Indeed, by H\"older inequality we have \\
  {\small$$
  \begin{aligned}
 \left[ u_{\varepsilon}-u\right]_{(s,M)}&=\sup_{\rho(v,\overline{M})\leqslant 1}\left|  \int_{\R^N}\int_{\R^N}\dfrac{(u_{\varepsilon}(x)-u(x))-(u_{\varepsilon}(y)-u(y))}{|x-y|^sM^{-1}(|x-y|^N)}v(x,y)dxdy\right|\\
  &\leqslant \sup_{\rho(v,\overline{M})\leqslant 1}\int_{\R^N} J(z)dz\int_{\R^N}\int_{\R^N}\left|  \dfrac{(u(x-\varepsilon z)-u(y- \varepsilon z))-(u(x)-u(y))}{|x-y|^sM^{-1}(|x-y|^N)}\right| \left| v(x,y)\right| dxdy\\
  & \leqslant \sup_{\rho(v,\overline{M})\leqslant 1} 2 ||v||_{{\overline{M}}}  \int_{|z|<1} J(z) \left| \left| \dfrac{(u(x-\varepsilon z)-u(y- \varepsilon z))-(u(x)-u(y))}{|x-y|^sM^{-1}(|x-y|^N)}\right| \right| _{M}dz,
  \end{aligned}
  $$}
  this implies that\\
  $$ \left[ u_{\varepsilon}-u\right]_{(s,M)}\leqslant 2\int_{|z|<1} J(z) \left| \left| \dfrac{(u(x-\varepsilon z)-u(y- \varepsilon z))-(u(x)-u(y))}{|x-y|^sM^{-1}(|x-y|^N)}\right| \right| _{M}dz.$$
  On the other hand, since $w(x,y):=\dfrac{|u(x)-u(y)|}{|x-y|^sM^{-1}(|x-y|^N)} \in E_M(\R^N \times \R^N )$, then given $\varepsilon>0$, there exists $g(x,y)\in C^{\infty}_0(\R^N \times \R^N)$ such that 
    $ ||w-g||_{M}\leqslant \dfrac{\varepsilon}{6}$. 
    That is, $$\left| \left| w(x-\epsilon z,y-\varepsilon z)-g(x-\varepsilon z,y-\varepsilon z)\right| \right| _{M}\leqslant \dfrac{\varepsilon}{6},$$ and for sufficiently small $\varepsilon$, $$\left| \left|g(x-\varepsilon z,y-\varepsilon z)-g(x,y)\right| \right|_{M}\leqslant\dfrac{\varepsilon}{6},$$ for every $z$ with $|z|\leqslant 1$. Thus $[u_{\varepsilon}-u]_{(s,M)}\leqslant \varepsilon$.\\
     \hspace*{15cm$\Box$ }
     
     Let $W^s_0{L_M}(\Omega)$
      denote the closure of  $C^{\infty}_0(\Omega)$ in the norm $||.||_{s,M}$ defined
     in $(\ref{6})$. The space $W^s_0{E_M}(\Omega)$ is defined in analogous fashion.\\ Note that, in view of Theorem \ref{3.2.},
     we have 
     $$W^s_0{E_M}(\R^N)=W^s{E_M}(\R^N).$$
    
       \begin{thm}(Generalized Poincar\'e inequality). 
         Let $\Omega$ be a bounded open subset of  $\R^N,$ and let $s\in (0,1)$. Let $M$ be an N-function. Then there exists a positive constant $\mu$ such that, \\         $$ ||u||_{M}\leqslant \mu [u]_{s,M}, \text{   } \forall u \in  W_0^sL_{M}(\Omega).$$
                        \end{thm}
          Therefore, if $\Omega$ is bounded and $M$ be an N-function, then $ [.] _ {s, M}$ is a norm of $W^s_{0}L_M (\Omega) $ equivalent to $ ||. ||_{s,M}.$
            \\

          \noindent \textbf{Proof}. 
          Since $W^s_0L_{M}(\Omega)$ is the closure of $C_0^{\infty}(\Omega)$ in $W^s{L_M}(\Omega)$, then it is enough to prove that there exists a positive constant $\mu$ such that, \\
          $$||u||_{M}\leqslant \mu [u]_{s,M}, \text{    } \forall u \in  C_0^{\infty}(\Omega).$$
          Indeed, let $ u \in C_0^{\infty}(\Omega) $ and $ B_R \subset \ R ^ N \smallsetminus \Omega $, i.e, the ball of radius $ R $ in the complement of $ \Omega $. Then for all $x\in \Omega$, $y\in B_R$ and all $\lambda>0$ we have, \\
          $$M(\dfrac{u(x)}{\lambda})=M\left( \dfrac{u(x)-u(y)}{\lambda |x-y|^s M^{-1}(|x-y|^N)}|x-y|^sM^{-1}(|x-y|^N)\right),$$
          this implies that,\\
          $$M(\dfrac{u(x)}{\lambda})\leqslant M\left( \dfrac{u(x)-u(y)}{\lambda |x-y|^sM^{-1}(|x-y|^N)} diam(\Omega\cup B_R)^s diam(M^{-1}(|\Omega\cup B_R|^N))\right),$$
we suppose $\alpha =diam(\Omega\cup B_R)^s diam(M^{-1}(|\Omega \cup B_R|^N))$, we get
$$M(\dfrac{u(x)}{\alpha\lambda})\leqslant M\left( \dfrac{u(x)-u(y)}{\lambda |x-y|^sM^{-1}(|x-y|^N)} \right),$$
therefore
         \small{ $$|B_R|M(\dfrac{u(x)}{\alpha\lambda})\leqslant \int_{B_R} M\left( \dfrac{u(x)-u(y)}{\lambda |x-y|^sM^{-1}(|x-y|^N)} \right) dy,$$}
then
           $$ \int_{\Omega} M(\dfrac{u(x)}{\alpha\lambda})dx \leqslant \dfrac{1}{|B_R|}\int_{\Omega}\int_{B_R} M\left( \dfrac{u(x)-u(y)}{\lambda |x-y|^sM^{-1}(|x-y|^N)} \right)dxdy,$$ 
           so, 
          $$ ||u||_{M}\leqslant \mu [u]_{s,M} \text{   } \forall u \in  C^{\infty}_0(\Omega),$$
          where $\mu =\dfrac{\alpha}{|B_R|}$. 
          By passing to the limit, the desired result is obtained.
          \\
               \hspace*{15cm$\Box$ }    
\subsection{Some embeddings results}$\label{101}$ 
The embeddings results obtained in the fractional Sobolev space $W^{s,p}(\Omega)$ can also be formulated for the fractional Orlicz-Sobolev spaces.

Let $M$ be a given N-function, satisfying the following conditions :
  \begin{equation}\label{15}
  \int_{0}^{1} \dfrac{M^{-1}(\tau)}{\tau^{\frac{N+s}{N}}}d\tau<\infty,
  \end{equation}
  
  \begin{equation}\label{16}
  \int_{1}^{\infty} \dfrac{M^{-1}(\tau)}{\tau^{\frac{N+s}{N}}}d\tau=\infty.
  \end{equation}
  For instance if $M(t)=\frac{1}{p}t^p$, then (\ref{15}) holds precisely when $sp<N$.\\
  If (\ref{16}) is satisfied, we define the inverse Sobolev conjugate N-function of $M$ as follows, 
  \begin{equation}\label{17}
  M_*^{-1}(t)=\int_{0}^{t}\dfrac{M^{-1}(\tau)}{\tau^{\frac{N+s}{N}}}d\tau.
  \end{equation}
\begin{thm}\label{3.4}
 Let $M$ be an $N$-function, and $s\in (0,1)$. Let $\Omega$  be a bounded open
  subset of  $\R^N$ with $C^{0,1}$-regularity 
    and bounded boundary. If $(\ref{15})$ and $(\ref{16})$  hold, then 
 \begin{equation}\label{18}
  W^s{L_M}(\Omega)\hookrightarrow L_{M_*}(\Omega).
 \end{equation}
\end{thm}
The proof will be carried out in a several lemmas. The first of these establishes
an estimate for the Sobolev conjugate N-function $M_*$, defined by (\ref{17}).
\begin{lem}$\label{3.2}$
  Let $M$ be an N-function satisfying (\ref{15}) and (\ref{16}), and suppose that,
  for some $p$ such that $1\leqslant p<N$,  the function $B$ defined by $B(t)=M(t^{\frac{1}{p}})$ is an N-function. Let $M_*$ be defined by (\ref{17}). Then the following
  conclusions may be drawn.

  \begin{enumerate}

  \item $[M_*(t)]^{\frac{N-s}{N}}$ is an N-function, in particular,  $M_*$ is an N-function.
  
  \item For every $\epsilon >0$, there exists a constant $ K_\epsilon>0$ such that for every $t$,
    \end{enumerate}
  \begin{equation}\label{19}
  [M_*(t)]^{\frac{N-s}{N}}\leqslant \dfrac{1}{2\epsilon}M_*(t)+\dfrac{K_\epsilon}{\epsilon}t.
  \end{equation}
   
  \end{lem}
 \noindent \textbf{Proof of lemma \ref{3.2}}. 
   (1) Let $Q(t)=[M_*(t)]^{\frac{N-s}{N}}$, Noting that $B^{-1}(t)=[M^{-1}(t)]^p$, we get 
   $$
   \begin{aligned}
   (Q^{-1})'(t)&=\dfrac{d}{dt} M^{-1}_*(t^{\frac{N}{N-s}}) \\
   &= \dfrac{N}{N-s} t^{\frac{N}{N-s}-1}\dfrac{ M^{-1}(t^{\frac{N}{N-s}})}{[t^{\frac{N}{N-s}}]^{\frac{N+s}{N}}}\\
    &= \dfrac{N}{N-s}\dfrac{M^{-1}(t^{\frac{N}{N-s}})}{t^{1+\frac{s}{N-s}}}\\
    &=\dfrac{N}{N-s}\left[ \dfrac{B^{-1}(t^{\frac{N}{N-s}})}{t^{\frac{N}{N-s}}} \right] ^{\frac{1}{p}} t^{-\mu},
   \end{aligned}
   $$
   where $\mu=1+\dfrac{s}{N-s}-\dfrac{N}{N-s}\dfrac{1}{p}=\dfrac{N(p-1)}{(N-s)p}\geqslant 0$. 
   Being
   the inverse of an N-function,   $B^{-1}$ satisfies $$\lim\limits_{t\rightarrow 0^+}\dfrac{B^{-1}(t)}{t} = \infty \text{ and } \lim\limits_{t\rightarrow \infty}\dfrac{B^{-1}(t)}{t} = 0,$$ and for $0<r<\sigma $ we have,  $\dfrac{B^{-1}(r)}{B^{-1}(\sigma)}>\dfrac{r}{\sigma}$. Hence, if  $0<t<s$, then we get,
   $$ \dfrac{(Q^{-1})'(t)}{(Q^{-1})'(s)}\geqslant (\dfrac{s}{t})^{-\mu}>1.$$
   It follows that $(Q^{-1})'$  is positive and decreases monotonically from $\infty$ to $0$ as t
   increases from $0$ to $\infty$, so that $Q$ is an N-function.\\
   (2)  Let $g(t)=\dfrac{M_*(t)}{t}$ and $h(t)=\dfrac{[M_*(t)]^{\frac{N-s}{N}}}{t}$. It readily checked that $h$ is bounded on finite intervals and $\lim\limits_{t\rightarrow \infty}\dfrac{g(t)}{h(t)}=\infty$, then for all $\varepsilon>0$  there exists $ t_0>0$ such that for every $t\geqslant t_0$, $h(t)\leqslant \dfrac{g(t)}{2\varepsilon}$. We pose $K_\varepsilon=\varepsilon \sup\limits_{0\leqslant t\leqslant t_0}h(t)$, then \\
   $$[M_*(t)]^{\frac{N-s}{N}}\leqslant \dfrac{1}{2\varepsilon}M_*(t)+\dfrac{K_\varepsilon}{\varepsilon}t.$$
     \hspace*{15cm $\Box$ }     
 \begin{lem}$\label{3.3}$
   Let $\Omega$ be an open subset of $\R^N$, and $0<s<1$. Let $f$ satisfies a
   Lipschitz-condition on $\R$ and $f(0)=0$, then,
   \begin{enumerate}
    \item For every $u\in W^{s,1}_{loc}(\Omega)$, $g\in W^{s,1}_{loc}(\Omega)$ where  $g(x)=f(|u(x)|)$.
    \item For every $u\in W^s{L_M}(\Omega)$, $g\in  W^s{L_M}(\Omega)$ where  $g(x)=f(|u(x)|)$.
       \end{enumerate}
   \end{lem}
    \noindent \textbf{Proof of lemma \ref{3.3}}. 
   (1) Let $K$ be a compact subset of $\Omega$, follows that $\mathds{1}_{K}g\in W^{s,1}(\Omega)$.\\
Since $f(0)=0,$ then we have, \\
$$\int_{\Omega}| \mathds{1}_K(x)g(x)|dx=\int_{\Omega} |\mathds{1}_K(x)(f(u(x))-f(0))|dx\leqslant C\int_{\Omega}|\mathds{1}_K(x)u(x)|dx<\infty ,$$
 where $C$ is the Lipschitz constant of $f$. On the other hand,
{\small\begin{equation}\label{20}
\begin{aligned}
    \int_{\Omega}\int_\Omega\left| \dfrac{\mathds{1}_K(x) g(x)-\mathds{1}_K(y) g(y)}{|x-y|^{N+s}}\right| dxdy = &\int_{K}\int_K\dfrac{ |g(x)-g(y)|}{|x-y|^{N+s}}dxdy +2\int_{\Omega \smallsetminus K}\int_{K}\dfrac{|\mathds{1}_K(x) g(x)|}{|x-y|^{N+s}}dxdy\\
    & +\int_{\Omega \smallsetminus K}\int_{\Omega \smallsetminus K}\dfrac{|\mathds{1}_K(x) g(x)-\mathds{1}_K(y) g(y)|}{|x-y|^{N+s}}dxdy,
    \end{aligned}
\end{equation}}

where the third term in the right hand-side of (\ref{20}) is null, and since $f$ satisfies the Lipschitz-condition, then,

   {\small \begin{equation}\label{21}
    \begin{aligned}
          \int_{\Omega}\int_\Omega\left| \dfrac{\mathds{1}_K(x) g(x)-\mathds{1}_K(y) g(y)}{|x-y|^{N+s}}\right| dxdy \leqslant &C\int_{K}\int_K\dfrac{| u(x)-u(y)|}{|x-y|^{N+s}}dxdy +2\int_{\Omega \smallsetminus K}\int_{ K}\dfrac{|g(x)|}{|x-y|^{N+s}}dxdy,\\
          \end{aligned}
    \end{equation} }
    where the first term in the right hand-side of (\ref{21}) is finite since
      $u\in W^{s,1}_{loc}(\Omega)$ and
      $$2\int_{K}\int_{\Omega\smallsetminus K}\dfrac{|g(x)|}{|x-y|^{N+s}}dxdy\leqslant 2\int_{K}|g(x)|dx\int_{\Omega\smallsetminus K}\dfrac{1}{d(y,\partial K)^{N+s}}dy< \infty.$$
     Note that due to the fact that $K$ is a compact subset, then $dis(y,\partial K)^{N+s}>0$  for all $y \in \R^N\smallsetminus K$ and we have $N+s>N$.\\

Therefore,
      $$  \int_{\Omega}\int_\Omega \left| \dfrac{\mathds{1}_K(x) g(x)-\mathds{1}_K(y) g(y)}{|x-y|^{N+s}}\right| dxdy< \infty.$$
      (2) Let $u\in W^s{L_M}(\Omega)$ then 
    there exists  $\lambda >0$ such that, \\
      $$ \int_{\Omega} \int_{\Omega} M(\dfrac{\lambda u(x)-\lambda u(y)}{|x-y|^sM^{-1}(|x-y|^N)})dxdy< \infty. $$
   Let $C>0$ denotes the Lipschitz constant of $f$ then, \\
   if $|C|\leqslant 1$  we have, 
   $$ \begin{aligned}
   \int_{\Omega} \int_{\Omega} M(\dfrac{\lambda g(x)-\lambda g(y)}{|x-y|^sM^{-1}(|x-y|^N)})dxdy &=\int_{\Omega} \int_{\Omega} M(\dfrac{\lambda f\circ u(x)-\lambda f\circ u(y)}{|x-y|^sM^{-1}(|x-y|^N)})dxdy\\
   &\leqslant \int_{\Omega} \int_{\Omega} M(\dfrac{|C|(\lambda u(x)-\lambda u(y))}{|x-y|^sM^{-1}(|x-y|^N)})dxdy\\
   & \leqslant \int_{\Omega} \int_{\Omega} M(\dfrac{\lambda u(x)-\lambda u(y)}{|x-y|^sM^{-1}(|x-y|^N)})dxdy<\infty.\\
   \end{aligned}
   $$
   If $|C|>1$ for $ \lambda_1=\dfrac{\lambda}{|C|}$ we get,
  {\small $$ \begin{aligned}
    \int_{\Omega} \int_{\Omega} M(\dfrac{\lambda_1 g(x)-\lambda_1 g(y)}{|x-y|^sM^{-1}(|x-y|^N)})dxdy &=\int_{\Omega} \int_{\Omega} M(\dfrac{\lambda_1 f\circ u(x)-\lambda_1 f\circ u(y)}{|x-y|^sM^{-1}(|x-y|^N)})dxdy\\
    &\leqslant \int_{\Omega} \int_{\Omega} M(\dfrac{|C|(\lambda u(x)-\lambda u(y))}{|C||x-y|^sM^{-1}(|x-y|^N)})dxdy\\
    & \leqslant \int_{\Omega} \int_{\Omega} M(\dfrac{\lambda u(x)-\lambda u(y)}{|x-y|^sM^{-1}(|x-y|^N)})dxdy<\infty,\\
    \end{aligned}
    $$}
    this implies that $g\in W^s{L_M}(\Omega)$.
      \\
      \hspace*{15cm$\Box$ } \\   
   \hspace*{0.3cm}  Let $M$ be an N-function, since $\lim\limits_{t\rightarrow 0}\frac{M(t)}{t}=0$ so there exists $\alpha>0$ such that $M(t)\leqslant t$ for all $t\leqslant \alpha$. For this $\alpha$, we define the function $M_1$ as, 
     
      \begin{equation}\label{22.}
      M_1(t)= \left \{
       \begin{array}{clclc}
      \frac{M(\alpha)}{\alpha} t\hspace{1cm}  & if & t\leqslant \alpha,\\\\
         M(t) \hspace*{1cm} & if & t>\alpha.
       \end{array}
       \right .  
       \end{equation}
   $M_1$ is a convex, continuous, nondecreasing, finite valued  function which is $M_1(0)=0$ and  \small{$\lim\limits_{t\rightarrow +\infty} M_1(t)=+\infty$}. $M_1$ is called a Young function (cf. \cite{33.} ).\\ 
    For a given domain $\Omega$ in $\R^N$, we define the space $L_{M_1}(\Omega)$ as,
    $$L_{M_1}(\Omega) =\Bigg\{u : \Omega\rightarrow \R : \exists\lambda>0  / \int_{\Omega}M_1( \lambda u(x))dx < \infty\Bigg\},$$ 
   this space is equipped with the norm,    
   
      \begin{equation}
       ||u||_{M_1}=\inf \Bigg\{\lambda>0  : \int_{\Omega}M_1( \dfrac{u(x)}{\lambda})dx \leqslant 1\Bigg\}.
      \end{equation}
   The Young complement of $M_1$ is defined for $0\leqslant x < \infty$ by $$\overline{M_1}(t)=\max \limits_{s\geqslant 0}\left\lbrace st-M_1(s) \right\rbrace, $$ then we have $st\leqslant M_1(t)+\overline{M_1}(s)$ for all $s,t\geqslant 0$ and for all $u,v\in L_{M_1}(\Omega)$ we get the H\"older inequality, i.e,
   $$\int_{\Omega}| u(x)v(x)|dx\leqslant 2||u||_{M_1}||v||_{\overline{M}_1} $$
\begin{lem}\label{3.3.}
Let $\Omega$ be a bounded open subset of $\R^N$ and let $s\in (0,1)$. Let $M$ be an N-function and $M_1$ as defined by $(\ref{22.})$ then,
\begin{enumerate}
\item $L_{M_1}(\Omega)=L_{M}(\Omega)$.
\item The norm $||.||_M$ and $||.||_{M_1}$ are equivalent.
\end{enumerate}
\end{lem}   
  \noindent \textbf{Proof of lemma \ref{3.3.}}. 
  (1) By definition of the function $M_1$ we have, $M(t)\leqslant\beta M_1(t)$ for all $t>0$, where $\beta=\max\left\lbrace 1, \frac{\alpha}{M(\alpha)}\right\rbrace$, then 
  $$L_{M_1}(\Omega)\subset L_M(\Omega).$$\\
  Let $u\in L_M(\Omega)$, we get 
  $$
                \begin{aligned}
 \int_{\Omega} M_1(\lambda u(x))dx&= \int_{\Omega\cap\left\lbrace \lambda u(x)\leqslant \alpha\right\rbrace } M_1(\lambda u(x))dx+ \int_{\Omega\cap \left\lbrace  \lambda u(x)>\alpha\right\rbrace } M_1(\lambda u(x))dx  \\
 &\leqslant M(\alpha)|\Omega|+  \int_{\Omega} M(\lambda u(x))dx< \infty.
 \end{aligned} $$            
 Then, $L_{M_1}(\Omega)=L_{M}(\Omega)$.\\ 
 (2) Let $u\in L_{M_1}(\Omega)$, since $M(t)\leqslant \beta M_1(t)$ for all $t>0$, then evidently $ ||u||_M\leqslant \beta||u||_{M_1}$.\\
 On the other hand, we get
 $$
   \begin{aligned}
  \int_{\Omega} M_1(\dfrac{u(x)}{||u||_M})dx&= \int_{\Omega\cap\left\lbrace \frac{u(x)}{||u||_M}\leqslant \alpha\right\rbrace } M_1(\dfrac{u(x)}{||u||_M})dx+ \int_{\Omega\cap \left\lbrace  \frac{u(x)}{||u||_M}>\alpha\right\rbrace } M_1(\dfrac{u(x)}{||u||_M})dx  \\
   &\leqslant M(\alpha) |\Omega|+  \int_{\Omega} M(\dfrac{u(x)}{||u||_M})dx< |\Omega|+1,              
                 \end{aligned}                
                           $$ 
so, $||u||_{M_1}\leqslant(M(\alpha)|\Omega|+1) ||u||_{M}$.\\
     \hspace*{15cm $\Box$ }

\begin{rem}
Let $\Omega$ be an open subset of $\R^N$ and let $s\in (0,1)$. Let $M_1$ as defined by $(\ref{22.})$, then we define the space $W^sL_{M_1}(\Omega)$ by,
\small{\begin{equation}
W^s{L_{M_1}}(\Omega)=\Bigg\{u\in L_{M_1}(\Omega) : \exists \lambda >0 / \int_{\Omega} \int_{\Omega} M_1\left( \dfrac{\lambda( u(x)- u(y))}{|x-y|^sM_1^{-1}(|x-y|^N)}\right) dxdy< \infty \Bigg\}.
\end{equation}
}
which equipped with the norm
$$||u||_{s,M_1}=||u||_{M_1}+[u]_{s,M_1}$$
where, $$[u]_{s,M_1}=\inf \Bigg\{\lambda >0 :  \int_{\Omega} \int_{\Omega} M_1\left( \dfrac{u(x)- u(y)}{\lambda|x-y|^sM_1^{-1}(|x-y|^N)}\right) dxdy\leqslant 1 \Bigg\}.$$
If  $\Omega$ is a bounded open subset of $\R^N$, then by lemma $\ref{3.3.}$, we have $W^sL_{M_1}(\Omega)=W^sL_{M}(\Omega)$ and the norm $||.||_{s,M}$ and $||.||_{s,M_1}$ are equivalent.
\end{rem}
\begin{lem}\label{3.4.}
Let $\Omega$ be a bounded open subset of $\R^N$ and let $s\in (0,1)$. Let $M_1$ as defined by $(\ref{22.})$, then the space $W^sL_{M_1}(\Omega)$ continuously embedded in $W^{s,1}(\Omega)$. Therefore $W^sL_{M}(\Omega)$ continuously embedded in $W^{s,1}(\Omega)$.
\end{lem}
 \noindent \textbf{Proof of lemma \ref{3.4.}}. Let $u\in W^sL_{M_1}(\Omega)$, we have by H\"older inequality
\begin{equation}\label{25.}
\int_{\Omega}|u(x)|dx\leqslant 2||u||_{M_1}||1||_{\overline{M_1}}.
\end{equation}
  On the other hand, we get
  \small{$$
     \begin{aligned}
  \int_{\Omega}\int_{\Omega} \dfrac{|u(x)-u(y)|}{|x-y|^{s+N}}dxdy&=  \int_{\Omega}\int_{\Omega\cap\left\lbrace  |x-y|\leqslant \alpha\right\rbrace } \dfrac{|u(x)-u(y)|}{|x-y|^{s+N}}dxdy+  \int_{\Omega}\int_{\Omega\cap\left\lbrace |x-y|> \alpha\right\rbrace } \dfrac{|u(x)-u(y)|}{|x-y|^{s+N}}dxdy\\
  &= I_1+I_2.
     \end{aligned}
$$}
By definition of $M_1$ and H\"older inequality, we have
\begin{equation}\label{26.}
     \begin{aligned}
  I_1=\int_{\Omega}\int_{\Omega\cap\left\lbrace  |x-y|\leqslant \alpha\right\rbrace } \dfrac{|u(x)-u(y)|}{|x-y|^{s+N}}dxdy&=\int_{\Omega}\int_{\Omega\cap\left\lbrace  |x-y|\leqslant \alpha\right\rbrace } \dfrac{|u(x)-u(y)|}{|x-y|^{s}M_1^{-1}(|x-y|^N)}dxdy\\
 & \leqslant 2[u]_{s,M_1}||1||_{\overline{M_1}},
     \end{aligned}
\end{equation}
and 
$$
     \begin{aligned}
  I_2&=\int_{\Omega}\int_{\Omega\cap\left\lbrace  |x-y|> \alpha\right\rbrace } \dfrac{|u(x)-u(y)|}{|x-y|^{s+N}}dxdy\\
  &=\int_{\Omega}\int_{\Omega\cap\left\lbrace  |x-y|> \alpha\right\rbrace } \dfrac{|u(x)-u(y)|}{|x-y|^{s}M_1^{-1}(|x-y|^N)}\dfrac{M_1^{-1}(|x-y|^N)}{|x-y|^N}dxdy\\
 & \leqslant \sup_{\Omega \times \Omega \cap\left\lbrace  |x-y|> \alpha\right\rbrace}\dfrac{M_1^{-1}(|x-y|^N)}{|x-y|^N}\int_{\Omega}\int_{\Omega\cap\left\lbrace  |x-y|> \alpha\right\rbrace } \dfrac{|u(x)-u(y)|}{|x-y|^{s}M_1^{-1}(|x-y|^N)}dxdy,
     \end{aligned}
$$
since $M^{-1}_1(t)$ is continuous for all $t>\alpha$ and $\Omega$ is bounded so,
$$\sup_{\Omega \times \Omega \cap\left\lbrace  |x-y|> \alpha\right\rbrace}\dfrac{M^{-1}(|x-y|^N)}{|x-y|^N}=C'< \infty.$$
Therefore by H\"older inequality, 
\begin{equation}\label{27.}
I_2\leqslant C'\int_{\Omega}\int_{\Omega\cap\left\lbrace  |x-y|> \alpha\right\rbrace } \dfrac{|u(x)-u(y)|}{|x-y|^{s}M_1^{-1}(|x-y|^N)}dxdy\leqslant 2C' [u]_{s,M_1}||1||_{\overline{M_1}}.
\end{equation}
Combining (\ref{25.}), (\ref{26.}) and (\ref{27.}) we obtain
$$||u||_{W^{s,1}}\leqslant C ||u||_{s,M_1},$$
where $C=(2+2C')||1||_{\overline{M_1}}$.\\
     \hspace*{15cm $\Box$ }     

            \noindent \textbf{Proof of theorem \ref{3.4}}.
              Let $\sigma(t)=[M_*(t)]^{\frac{N-s}{N}}$
              and $u\in W^s{L_M}(\Omega)$, we suppose for the moment that $u$ is bounded on $\Omega$ and not equal to zero in $L_M(\Omega)$, then  $\displaystyle\int_\Omega M_*(\dfrac{u(x)}{\lambda})dx$  decreases continuously from infinity to zero as  $\lambda$ increases from zero to infinity, so that 
        \begin{equation}\label{22}
         \int_\Omega M_*(\dfrac{u(x)}{k})dx =1 \text{ , } k=||u||_{M_*}.
        \end{equation}
              Let $f(x)=\sigma(\dfrac{u(x)}{k})$. Evidently by lemma \ref{3.4.} $u\in W^{s,1}(\Omega)$, and $\sigma$ is Lipschitz, so that by lemma \ref{3.3} we have $f\in W^{s,1}(\Omega)$, and  since $N>s$, then by theorem \ref{2.3}, one has, \\
              $$W^{s,1}(\Omega) \hookrightarrow L^{\frac{N}{N-s}}(\Omega).$$ 
              So 
              $$||f||_{L^{\frac{N}{N-s}}}\leqslant k_1 \left( ||f||_{L^1}+[f]_{s,1}\right), $$
              and by (\ref{22}), 
              $$1= \left( \int_\Omega M_*(\dfrac{u(x)}{k})dx \right) ^{\frac{N-s}{N}}=||f||_{L^{\frac{N}{N-s}}},$$
              this implies that, \\
             \begin{equation}\label{23}
              \begin{aligned}
                 1 &\leqslant  k_1 \left( ||f||_{L^1}+[f]_{s,1}\right)\\
                 &= k_1 \left( \int_{\Omega}\sigma(\dfrac{u(x)}{k})dx+\int_{\Omega}\int_\Omega\dfrac{ |f(x)-f(y)|}{|x-y|^{N+s}}dxdy\right) \\
      &=  k_1 \left( \int_{\Omega}\sigma(\dfrac{u(x)}{k})dx+\int_{\Omega}\int_\Omega\dfrac{|\sigma(\dfrac{u(x)}{k})-\sigma(\dfrac{u(y)}{k})|}{|x-y|^{N+s}}dxdy\right) \\
                & = k_1 I_1 +k_1 I_2.
                   \end{aligned}   
             \end{equation}
            By (\ref{19}) we have for $\varepsilon =k_1$,
        \begin{equation}\label{24}
        k_1I_1 \leqslant \dfrac{1}{2} \int_\Omega M_*(\dfrac{u(x)}{k})dx+\dfrac{k_\varepsilon}{k}\int_\Omega |u(x)|dx \leqslant \dfrac{1}{2} +\dfrac{k'_\varepsilon}{k}||u||_M,
        \end{equation}
         where $k'_\varepsilon =2k_\varepsilon ||1||_{\overline{M}}$ since $\Omega$ has a finite volume.\\
          On the other hand, since $\sigma$ is Lipschitz, then there exists $C>0$ such that, \\
         $$ k_1I_2\leqslant \dfrac{C}{k}\int_{\Omega}\int_\Omega\dfrac{ |u(x)-u(y)|}{|x-y|^{N+s}}dxdy. $$
         But by the lemma \ref{3.4.}, we have
        \begin{equation}\label{25}
        \int_{\Omega}\int_\Omega\dfrac{ |u(x)-u(y)|}{|x-y|^{N+s}}dxdy\leqslant C'[u]_{s,M},
        \end{equation}
        and \\
          \begin{equation}\label{26}
          k_1I_2\leqslant \dfrac{C}{k} C'[u]_{s,M}.
          \end{equation}
          We pose $k_3=Ck_1 C'$. Combining (\ref{24})-(\ref{26}) we obtain \\
          $$1\leqslant  \dfrac{1}{2} +\dfrac{k'_\varepsilon}{k}||u||_M+\dfrac{k_3}{k}[u]_{s,M},$$
          this implies that,
          $$ \dfrac{k}{2} \leqslant k'_\varepsilon ||u||_M + k_3[u]_{s,M}.$$
         So we obtain, \\
          $$||u||_{M_*}\leqslant k_4||u||_{s,M},$$
          where $k_4= \max \left\lbrace  2k'_{\varepsilon} , 2 k_3\right\rbrace $.\\
           If $u\in W^s{L_M}(\Omega)$ arbitrary, we define  
           
           $$ u_n(x)=\left \{
           \begin{array}{lcl}
          u(x)\hspace*{1cm} \text{ if } |u(x)|\leqslant n,\\
           n \text{ sgn } u(x) \text{ if }|u(x)|>n.\\
           \end{array}
           \right .$$
           
           $u_n$ is bounded
     and  by lemma 3.4 it is belongs to  $W^s{L_M}(\Omega)$. Moreover $$||u_n||_{M_*}\leqslant k_4||u_n||_{s,M}\leqslant k_4||u||_{s,M}.$$
           \\ Let $\lim\limits_{n \rightarrow \infty}||u_n||_{M_*}=k$, then $ k\leqslant k_4||u||_{s,M}$. By Fatou's Lemma we get \\
           $$ \int_\Omega M_*(\dfrac{u(x)}{k})dx \leqslant \lim\limits_{n\rightarrow \infty}\int_\Omega M_*(\dfrac{u_n(x)}{||u_n||_{M_*}})dx<1, $$
           so $u\in L_{M_*}(\Omega)$ and $||u||_{M_*}\leqslant k.$
\\            
  \hspace*{15cm$\Box$ }     
          
         \begin{thm}
        Let $s\in (0,1)$ and $M$ be an $N$-function. Let $\Omega$  be a bounded open
          subset of  $\R^N$ and  $C^{0,1}$-regularity 
            with bounded boundary. If $(\ref{15})$ and  $(\ref{16})$ hold, then 
         \begin{equation}\label{27}
          W^s{L_M}(\Omega)\hookrightarrow L_{B}(\Omega),
         \end{equation}
         is compact for all $B\prec\prec M_*$.
         \end{thm}
          \noindent \textbf{Proof}.  
        By the lemma \ref{3.4.}, we have,
        $$ W^s{L_M}(\Omega)\hookrightarrow W^{s,1}(\Omega)\hookrightarrow L^1(\Omega).$$
       The latter embedding being compact by theorem \ref{2.4}. A bounded subset $S$ of $W^s{L_M}(\Omega)$ is also a bounded subset of $L_{M_*}(\Omega)$ and precompact in $L^1(\Omega)$, hence  by theorem  \ref{2.1,}  it is precompact in $L_B(\Omega)$.\\
  \hspace*{15cm$\Box$ }     
     \subsection{Application}
    In this final subsection, we define the fractional M-Laplacian operator, and we establish the existence of a unique solution for the variational problem related to this operator by  the Minty Browder theorem.\\
    
    In the rest of this subsection we assume that $(M,\Omega)$ is $\Delta$-regular.
    \begin{defini}
    Let $M$ be an N-function and $0<s<1$, we define the fractional M-Laplacian operator as, 
       {\small $$
        \begin{aligned}
        (-\Delta)^s_mu(x)&=2P.V \int_{\R^N} M'\left( \dfrac{(u(x)-u(y))}{|x-y|^sM^{-1}(|x-y|^{N}) }\right) \dfrac{u(x)-u(y)}{|u(x)-u(y)|}\dfrac{dy}{|x-y|^s M^{-1}(|x-y|^{N})}\\
        &=2P.V \int_{\R^N} m\left( \dfrac{(u(x)-u(y))}{|x-y|^sM^{-1}(|x-y|^{N}) }\right) \dfrac{u(x)-u(y)}{|u(x)-u(y)|}\dfrac{dy}{|x-y|^sM^{-1}(|x-y|^{N}) }
        \end{aligned}
         $$}
         where P.V is the principal value and $M'=m$.\\
        
    \end{defini}
 In the case $M(t)=\dfrac{|t|^p}{p}$, we have 
$$ (-\Delta)^s_mu(x)=(-\Delta)^s_p u(x)=2PV\int_{\R^N} \dfrac{|u(x)-u(y)|^{p-2}(u(x)-u(y))}{|x-y|^{N+sp}}dy .$$ 
\begin{lem}\label{3.5.}
If $u\in W^sL_M(\R^N)$, then $ (-\Delta)^s_mu(x) \in (W^sL_M(\R^N))^*$,
 
   and 
    $$<(-\Delta)^s_mu,v> = \int_{\R^N} \int_{\R^N} m(h_{x,y}(u))\dfrac{u(x)-u(y)}{|u(x)-u(y)|}h_{x,y}(v)dxdy, $$
        
 for all  $v\in W^sL_M(\R^N)$, where  $h_{x,y}(u):=\dfrac{(u(x)-u(y))}{|x-y|^sM^{-1}(|x-y|^{N}) }$.
\end{lem}
          \noindent \textbf{Proof}. 
 Evidently 
 $$<(-\Delta)^s_mu,v> = P.V\int_{\R^N} \int_{\R^N} m(h_{x,y}(u))\dfrac{u(x)-u(y)}{|u(x)-u(y)|}h_{x,y}(v)dxdy. $$
   By Young inequality and (\ref{2}), we get
  $$
  \begin{aligned}
  \left|  m(h_{x,y}(u))\dfrac{u(x)-u(y)}{|u(x)-u(y)|}h_{x,y}(v)\right| 
   &\leqslant \left| m(h_{x,y}(u))h_{x,y}(v)\right| \\
  & \leqslant \overline{M}\left(  m(h_{x,y}(u))\right) + M\left( h_{x,y}(v)\right) \\
  & \leqslant (p-1) M\left( h_{x,y}(u)\right) + M\left(h_{x,y}(v)\right) \in L^1(\R^N\times \R^N).
  \end{aligned}
  $$ 
  Finally 
   $$<(-\Delta)^s_mu,v> = \int_{\R^N} \int_{\R^N} m(h_{x,y}(u))\dfrac{u(x)-u(y)}{|u(x)-u(y)|}h_{x,y}(v)dxdy. $$ 
     \hspace*{15cm$\Box$ }\\     
  Given a bounded open set $\Omega \subset \R^N$, we establish the existence of unique weak solution for the following Dirichlet type equation,
  
  \begin{equation}\label{28}
  \left \{
  \begin{array}{clclc}
  (-\Delta)^s_mu&=&\hspace{-1cm} f & in & \Omega,\\\\
   \hspace{1cm}u&=& 0 \hspace*{1cm} & in & \R^N\smallsetminus \Omega.
  \end{array}
  \right .  
  \end{equation}
  We
  shall work in the closed linear subspace 
  $$\widetilde{W}^s_0L_M(\Omega)=\left\lbrace u\in W^sL_M(\Omega) \text{ : } u=0 \text{ a.e in } \R^N\smallsetminus \Omega \right\rbrace, $$
  equivalently renormed by setting $[.]_{s,M}$. which is a reflexive separable Banach space.
  \begin{defini}
  We say that $u\in \widetilde{W}^s_0L_M(\Omega)$ is a weak solution of (\ref{28}) if 
  \begin{equation}\label{29}
<(-\Delta)^s_mu,v> =\int_{\Omega}fvdx
  \end{equation}
  for all $v \in \widetilde{W}^s_0L_M(\Omega)$.
  \end{defini}
  \begin{thm}
  Let $\Omega$ be a bounded open subset of $\R^N$ and $f\in (\widetilde{W}^s_0L_M(\Omega))^*$, then the problem  (\ref{28}) has a unique solution $u\in \widetilde{W}^s_0L_M(\Omega)$.

  \end{thm}
            \noindent \textbf{Proof}. 
 We need to show that $(-\Delta)^s_m$ satisfies the conditions of Minty Browder theorem (cf. \cite{6}).\\
 \noindent \textbf{Step 1}. $(-\Delta)^s_m$ is bounded and continuous. 
 Indeed, by lemma $\ref{3.5.}$ and H\"older inequality we have for all $u \in \widetilde{W}^s_0L_M(\Omega)$
 $$||(-\Delta)^s_m u||_{(\widetilde{W}^s_0L_M(\Omega))^*}=\sup_{||v||_{s,M}\leqslant 1}<(-\Delta)^s_mu,v>\leqslant 2||m(h_{x,y}(u))||_{\overline{M}}$$
 therefore  $||(-\Delta)^s_m u||_{(\widetilde{W}^s_0L_M(\Omega))^*}$ is bounded once $||u||_{s,M}$ is bounded.\\
 \hspace*{0.3cm}Let $u_n \rightarrow u$ in $W^s_0L_M(\Omega)$ we show that $(-\Delta)^s_mu_n \rightarrow (-\Delta)^s_mu$ in $(\widetilde{W}^s_0L_M(\Omega))^*$. Indeed H\"older  inequality
{\small$$
    \begin{aligned}    
 ||(-\Delta)^s_mu_n-(-\Delta)^s_mu||_{(\widetilde{W}^s_0L_M(\Omega))^*}\leqslant ||m(h_{x,y}(u_n))\dfrac{u_n(x)-u_n(y)}{|u_n(x)-u_n(y)|}-m(h_{x,y}(u))\dfrac{u(x)-u(y)}{|u(x)-u(y)|} ||_{\overline{M}},
        \end{aligned} $$}
On the other hand, since $u_n \rightarrow u$ in $L_M(\Omega)$, by dominated convergence theorem, there exists a subsequence $\left\lbrace u_{n_k}\right\rbrace $ and a function $h$ in $L_M(\Omega)$ such that 
 $u_{n_k}(x)\rightarrow u(x)$ and 
 $|u_{n_k}(x)|\leqslant |h(x)|$ for all $k$, a.e. on $\Omega$. This implies that
 $$|m(h_{x,y}(u_{n_k}))\dfrac{u_{n_k}(x)-u_{n_k}(y)}{|u_{n_k}(x)-u_{n_k}(y)|}|\leqslant |m(h_{x,y}(h))| \in L_{\overline{M}}(\Omega) \text{ a.e in  } \Omega \times \Omega,$$
 and
 $$m(h_{x,y}(u_{n_k}))\dfrac{u_{n_k}(x)-u_{n_k}(y)}{|u_{n_k}(x)-u_{n_k}(y)|}\longrightarrow m(h_{x,y}(h)) \dfrac{u(x)-u(y)}{|u(x)-u(y)|}  \text{ a.e in } \Omega \times \Omega,$$ 
 then by dominated convergence theorem we obtain the desired result.\\
\noindent \textbf{Step 2}. $(-\Delta)^s_m$ is strictly monotonous. Since $m$ is increasing, then $f(u):=m(u)\dfrac{u}{|u|}$ is also increasing, then for all $u,v \in \widetilde{W}^s_0L_M(\Omega)$ such that $u \neq v$, we have 
{\small $$<(-\Delta)^s_mu -(-\Delta)^s_mv, u-v>=\int_{\Omega}\int_{\Omega}\left[  m(h_{x,y}(u))k_{x,y}(u)-m(h_{x,y}(v))k_{x,y}(v)\right] \left( h_{x,y}(u)-h_{x,y}(v)\right)> 0, 
$$}
where $k_{x,y}(u):=\dfrac{u(x)-u(y)}{|u(x)-u(y)|}$.\\
\noindent \textbf{Step 3}. $(-\Delta)^s_m$ is coercive. Indeed, let $\beta \in (1,[u]_{s,M})$,  By lemma C.3(ii) in \cite{16} we have 
$$
 \begin{aligned} 
\phi(u)&=\int_{\Omega} \int_{\Omega} M\left( \dfrac{\lambda u(x)-\lambda u(y)}{|x-y|^sM^{-1}(|x-y|^N)}\right) dxdy\\
&\geqslant \beta^{p0}\int_{\Omega} \int_{\Omega} M\left( \dfrac{\lambda u(x)-\lambda u(y)}{\beta|x-y|^sM^{-1}(|x-y|^N)}\right) dxdy\geqslant \beta^{p0}.
 \end{aligned} 
$$
On the other hand, we have $<\phi'(u),v>=<(-\Delta)^s_mu,v>$ and  since $M$ is convex, it follows that $\phi$ is also convex. Thus, we have
$$\phi(u)\leqslant <\phi'(u),u> \text{ for all } \in W^s_0L_M(\Omega),$$
it is clear that for any $u\in \widetilde{W}^s_0L_M(\Omega)$ with $[u]_{s,M}>1$ we have
$$\dfrac{<(-\Delta)^s_mu,v>}{||u||_{s,M}}=\dfrac{<\phi'(u),u>}{||u||_{s,M}}\geqslant \dfrac{\phi(u)}{||u||_{s,M}}\geqslant\dfrac{[u]^{p0}_{s,M}}{||u||_{s,M}}\geqslant C||u||_{s,M}^{p0-1},$$
where $C$ it is the constant of the Poincar\'e-inequality.
Thus,
$$\lim\limits_{||u||_{s,M}\rightarrow \infty}\dfrac{<(-\Delta)^s_mu,v>}{||u||_{s,M}}=\infty,$$
i.e. $(-\Delta)^s_m$ is coercive.\\
 Hence, in light of Minty-Browder theorem then there exists a unique solution $u\in \widetilde{W}^s_0L_M(\Omega)$ of the problem  (\ref{28}).\\
     \hspace*{15cm$\Box$ }\\

\end{document}